\documentclass[12pt]{amsart}
\usepackage{graphicx}
\usepackage[english]{babel}
\usepackage{blindtext}
\usepackage[utf8]{inputenc}
\usepackage{fancyhdr}
\usepackage{amsmath}
\usepackage{amssymb}
\usepackage{mathtools}
\usepackage{tikz-cd}
\usepackage{enumitem}
\usepackage{latexsym,amsfonts,amssymb,amsthm,amsmath}
\usepackage{amssymb}

\def\C{\mathbb C}

\def\R{\mathbb R}
\def\Q{\mathbb Q}
\def\P{\mathbb P}
\def\O{\mathcal O}

\def\Z{\mathbb Z}

\def\ra{\rightarrow}

\newtheorem{thm}{Theorem}[section]

\newtheorem{cor}{Corollary}[section]
\newtheorem{defn}{Definition}[section]
\newtheorem{rem}{Remark}[section]
\newtheorem{prop}{Proposition}[section]

\begin{document}

\title{On units with complex Galois conjugates of equal absolute value}

\renewcommand\rightmark{On units with complex Galois conjugates of equal absolute value}
\renewcommand\leftmark{\c Stefan Deaconu}

\author{ \c Stefan Deaconu}

\begin{abstract} We prove, using a theorem of Northcott--Weil, that if a number field $K$ with $s\geq 1$ real embeddings and $2t\geq 2$ complex ones has a group of units $U\subset \O_K^*$ such that all elements in $U$ have all its  complex conjugates of same absolute value, then one necessarily has $t=1.$ This fact has an interesting implication in complex hermitian geometry, namely  it describes all Oeljeklaus--Toma manifolds carrying locally conformally K\"ahler structures. It is then shown that the same method works for the situation concerning the existence of pluriclosed metrics on Oeljeklaus--Toma manifolds.
\end{abstract}
\maketitle
\section{Introduction} In view of the Kronecker's unit theorem, it is a natural question to look at algebraic units whose complex Galois conjugates have equal absolute value. For instance, $\sqrt{\sqrt[3]{2}-1}$ has this property, yet not being a root of unity. Obviously, for a given number field $K$, the set of units of $\mathcal{O}_K^*$ with this property form a subgroup $U;$ letting $s$ to be the number of real embeddings of $K$ and $2t$ the number of complex embeddings (for simplicity, call the pair $(s,t)$ {\em the signature} of $K$), then the Dirichlet logarithmic embedding of units implies at once that the rank of $U$ must be at most $s.$

This leads to the following question:

\hfill

{\em Give examples of number fields $K$ with signature $(s, t)$ (with $s, t>0$) having subgroups $U$ as above of (maximal) rank $s.$} 

\hfill

The question arose also from a problem in complex geometry: namely, such number fields would have given examples of compact complex manifolds (the so called  "OT manifolds"), with interesting geometric structures, namely "l.c.K metrics" (see Section \ref{geom} for the precise definitions).

Some geometrical facts hinted that such examples of number fields must necessary have $t=1$, but the proof of this came along a series of papers. First it was proven in \cite{Vu} that one must have $t<s;$ next, this result was widely extended in \cite{Dub} where it was shown that the signature of any  number field containing a unit whose complex conjugates are of same absolute value must obey a relation of the form $s = (2t + 2m)q-2t$, for some integers $m\geq 0$ and $q\geq 2.$ This left some cases $(s, t)$ still open. Eventually, the remaining cases  were settled in \cite{DV}, hence proving that any field $K$ as in the above question must indeed have $t=1$.

The goal of the present paper is to give a simple, unified proof of the above results. Moreover, we show that the same techinque can be used in order to treat the analogous problem raised by searching for another class of metrics on Oeljeklaus--Toma manifolds, namely pluriclosed (or S.K.T) metrics. 

\section{Preliminaries}

\subsection{Heights on number fields. Congruence subgroups}

Given a number field $K$ having signature $(s, t)$ with $s\not=0\not=t$, label $\sigma_1,\dots, \sigma_s$ its real embeddings and $\sigma_{s+1},\dots, \sigma_{s+2t}$ its complex embeddings, where
\begin{align*}
\sigma_{s+t+k}=\overline{\sigma_{s+k}}
\end{align*}
 for all $k=1,\dots, t$. Recal the following: 

\begin{thm}\label{deginf}\cite[Appendix 2. Theorem 1]{M} For an extension of number fields $K/k$, any embedding of $k$ extends to exactly $[K:k]$ embeddings of $K$.
\end{thm}

Let $\mathcal{O}_K$ be the associated ring of algebraic integers and $\mathcal{O}_K^*$ be its group of units.

Next, let $\mathbb{V}(K)$ a complete set of inequivalent (i.e. not inducing the same topology) absolute values on $K$ (so both, archimedean and nonarchimedean). Recall that this can be described as follows. First, $\mathbb{V}(K)=\mathbb{V}_0(K)\cup\mathbb{V}_{\infty}(K)$, where $\mathbb{V}_0(K)$ is the set of nonarchimedean absolute values and $\mathbb{V}_{\infty}(K)$ that of archi\-medean ones. For each prime $\mathfrak{p}$ of $K$, this is lying above some prime $p\in\Q$, and one takes
\begin{align*}
|\text{ }|_{\mathfrak{p}}:K\rightarrow\R_+, |x|_{\mathfrak{p}}=p^{\frac{-\text{ord}_{\mathfrak{p}}(x)}{e_{\mathfrak{p}}}},
\end{align*}
where $e_{\mathfrak{p}}=\text{ord}_{\mathfrak{p}}(p)$. Of course, the order of $\mathfrak{p}$ in $a\in\O_K$ means the exponent with which $\mathfrak{p}$ appear in the decomposition of $\langle a\rangle$ into prime factors (\cite[Chapter 3.Theorem 16]{M}), and then it is extended to any $x=\frac{a}{b}\in K$ by taking the difference, as usual. Thus, any prime ideal gives a nonarchimedean absolute value. Also, any embedding $\sigma:K\hookrightarrow\C$ is establishing an archimedean absolute value,
\begin{align*}
|\text{ }|_{\sigma}:K\rightarrow\R_+, |x|_{\sigma}=|\sigma(x)|,
\end{align*}
where $|\text{ }|$ mean the usual euclidean absolute value on $\R$ or $\C$. Modulo equivalence, Ostrowski's theorem says that these two ways produce all the absolute values on $K$. \\

\begin{rem}\label{cont} The absolute values at infinity are continuous. More precisely, let $|\text{ }|_{\sigma}\in\mathbb{V}_{\infty}(K)$ and let $K_{\sigma}$ be the completion of $K$ with respect to $|\text{ }|_{\sigma}$. Then there is a uniquely determined extension of $|\text{ }|_{\sigma}$ to a continuous absolute value $|\text{ }|_{\bar{\sigma}}: K_{\sigma}\rightarrow\R_+$; the field $K_{\sigma}$ is topologically isomorphic to $\R$ if $\sigma$ is a real embedding, or to $\C$ if $\sigma$ is a complex one \cite[1.2.3]{AA}. 
\end{rem}

Now let $P=[x_0:\ldots :x_n]\in\P^n(K)$.  {\em The height} of the point $P$ is defined as the quantity
\begin{equation}\label{height}
H_K(P)=\prod_{\nu\in\mathbb{V}(K)}\text{max}\{||x_0||_{\nu},\ldots ,||x_n||_{\nu}\},
\end{equation}
where $||x||_{\nu}$ is the normalized absolute value $|x|_{\nu}^{n_{\nu}}$, $n_{\nu}=[K_{\nu}:\Q_{\nu}]$ being the degree of the corresponding extension of completions (by $\nu$ it is meant an equivalence class of absolute values, called a {\em  place}). The product formula shows that the height is well--defined indeed \cite[B.2.1]{S}.

The most important fact for the present note is the following: \\

\begin{thm}\label{NW} (Northcott--Weil, \cite[B.2.3]{S}). For any $C\geq 0$, the set
$$\{P\in\P^n(K); H_K(P)\leq C\}$$
is finite. 
\end{thm}

Finally, recall:

\begin{defn}\label{cong} Let $K$ be a number field. A subgroup $U$ of $\mathcal{O}_K^*$ is called a congruence subgroup if there is an ideal $\mathfrak{a}$ of $\mathcal{O}_K$,  such that $U(\mathfrak{a})=\{\varepsilon\in\mathcal{O}_K^*; \varepsilon\equiv 1 (\text{mod }\mathfrak{a})\}\subseteq U$. 
\end{defn}

Note that any congruence subgroup is of finite index in $\mathcal{O}_K^*$: this is because $U(\mathfrak{a})$ is the kernel of $\O_K^*\rightarrow (\mathcal{O}_K/\mathfrak{a}, \cdot), \varepsilon\mapsto\text{cls}(\varepsilon)$, and the quotient ring $\O_K/\mathfrak{a}$ is finite (as $\mathcal{O}_K$ is a number ring). The interesting fact is that the converse also holds:

\begin{thm}\label{chevy}(Chevalley, \cite{C}) A subgroup $U\subseteq\mathcal{O}_K^*$ is a congruence subgroup if and only if it has finite index.
\end{thm}

\subsection{Oeljeklaus--Toma manifolds}\label{geom}

These manifolds were introduced by K.Oeljeklaus and M. Toma in \cite{OeTo}, as a generalization to higher dimensions of the Inoue surfaces $S_M$ (\cite{In}). Very briefly, their construction goes as follows. Fix a number field $K$ of signature $(s, t)$, and label its embeddings such that $\sigma_1,\dots, \sigma_s$ are the real ones, while $\sigma_{s+t+i}=\overline{\sigma}_{s+i}$ for any $i=1,\dots, t.$ Let $\mathcal{O}_K$ be the ring of integers of $K,$ let $\mathcal{O}_K^*$ be the group of units of $\O_K$ and $\mathcal{O}_K^{*, +}$ the subgroup of $\mathcal{O}_K^*$ of {\em totally positive} units, that is elements $u\in \mathcal{O}_K^*$ such that $\sigma_i(u)>0$ for all $i=1,\dots, s.$ Letting ${\mathbb H}:=\{z\in \mathbb{C}\vert {\rm Im}(z)>0\}$, we see there are natural actions of $\mathcal{O}_K$ and respectively of $\mathcal{O}_K^{*,+}$ on ${\mathbb H}^s\times \mathbb{C}^t\subset \mathbb{C}^{s+t}$ by
$$a\cdot(z_1,\ldots ,z_{s+t}))\mapsto (z_1+\sigma_1(a),\ldots ,z_{s+t}+\sigma_{s+t}(a)).$$
and respectively
$$u\cdot(z_1,\ldots ,z_{s+t})\mapsto (\sigma_1(u)z_1,\ldots ,\sigma_{s+t}(u)z_{s+t}).$$
The combined resulting action of $\mathcal{O}_K^{*, +}\ltimes \mathcal{O}_K$ is however not discrete in general. Still, in \cite{OeTo} it is shown that one can always find subgroups $U\subset \mathcal{O}_K^{*, +}$ such that the action of $U\ltimes \mathcal{O}_K$ is discrete and cocompact: the resulting compact complex manifold is usually denoted $X(K, U)$ and is called an {\em Olejeklaus-Toma manifold} (OT, for short). 

As these manifolds do not admit K\" ahler metrics \cite[2.5]{OeTo}, it is natural to ask whether other natural metrics (do) exist on them. One of the most interesting candidates are the {\em locally conformally K\"{a}hler} metrics (l.c.K, for short): these are those whose associated $(1,1)-$ forms $\omega$ have the property
\begin{equation}\label{lck}
d\omega=\theta\wedge \omega
\end{equation}
for some closed $1-$form $\theta$ (for more details see \cite{DO}). Another candidates are the {\em strongly K\"{a}hler with torsion} metrics (S.K.T), or {\em pluriclosed}, i.e. those for which  $\partial\overline{\partial}\omega=0$. The existence of such metrics on OT manifolds $X(K, U)$ can be read off the Galois properties of the group of units $U$. More precisely, it was shown that: 

\begin{prop}\label{Bat} (\cite{Dub}, appendix by L. Battisti, and \cite[6.5]{OI})
An Oeljeklaus-Toma manifold $X(K, U)$ admits an l.c.K metric if and only if for any unit $u\in U$ one has
\begin{equation}\label{rest}
\vert \sigma_{s+1}(u)\vert=\dots=\vert \sigma_{s+t}(u)\vert
\end{equation} 
\end{prop}

\begin{prop}\label{oti}(\cite{O}) An Oeljeklaus-Toma manifold $X(K, U)$, for which $s\geq 1$, admits a pluriclosed metric if and only if 
\begin{equation}\label{eqval2}
	\begin{split}
&s\leq t,\\
& \sigma_i(u)\vert \sigma_{s+i}(u)\vert^2=1\text{ for }i\in\{1,\ldots ,s\},\\
& \vert \sigma_{s+i}(u)\vert^2=1\text{ for }i\in\{s+1,\ldots, t\},
\end{split}
\end{equation}
\end{prop}

The analogue quantitative restriction imposed by the existence of pluriclosed metrics turns out to be $s=t$ (this and other metrics are considered in \cite{metric}).

\section{The main results}. 

\begin{thm}\label{main}
 Let $K$ be a number field having signature $(s, t)$, with $s\geq 1$ and $t\geq 1$. If $U$ is a rank $s$ subgroup of $\mathcal{O}_K^*$ such that for each $u\in U$ one has $\left |\sigma_{s+i}(u)\right |=|\sigma_{s+j}(u)|$, for all $i, j$, then $t=1$. 
\end{thm}

\begin{proof}
Assume $t>1$. Denote by  $N:=K^{\text{nor}}$ the  normal closure of $K$ identified with  a fixed  embedding of it into $\mathbb{C}$ and by 
$$G_U:=\{u\in U\vert \sigma_{s+1}(u)=\dots =\sigma_{s+t}(u)\}.$$ We consider the map $U\ra \P^{t-1}(N)$ given by
\begin{equation}\label{fin}
 u\mapsto P_u=[\sigma_{s+1}(u):\ldots :\sigma_{s+t}(u)]. 
\end{equation}
Plainly, this induces an injection $\rho:U/G_U\hookrightarrow \P^{t-1}(N)$. We infer that its image is finite. To prove this, consider the height $H_N$ (see (\ref{height})). For each $u\in U$, we have
\begin{align*}
H_N(P_u)&=H_N\left(\left[1:\frac{\sigma_{s+2}(u)}{\sigma_{s+1}(u)}:\ldots :\frac{\sigma_{s+t}(u)}{\sigma_{s+1}(u)}\right]\right).
\end{align*}
As all $\frac{\sigma_{s+k}(u)}{\sigma_{s+1}(u)} \; (k\geq 2)$ are units, we get that 
$$\left |\left |\frac{\sigma_{s+k}(u)}{\sigma_{s+1}(u)}\right |\right |_{\nu}=1$$ for each $\nu\in\mathbb{V}_0(K)$; this implies 

\begin{equation}\label{infi}
H_N(P_u)=\prod_{\nu\in\mathbb{V}_\infty(N)}\text{max}\left\{1, \left |\left |\frac{\sigma_{s+2}(u)}{\sigma_{s+1}(u)}\right |\right |_{\nu},\ldots ,\left |\left |\frac{\sigma_{s+t}(u)}{\sigma_{s+1}(u)}\right |\right |_{\nu}\right\}.
\end{equation}

Now, for any infinite place $\nu$, since the completion $\left(N_{\nu}, \vert\vert \cdot\vert\vert_{\nu}\right)$ is topologically equivalent to $(\C, \vert \cdot\vert)$ ($t>1$), we see that, by letting $S$ to be the preimage in $N_{\nu}$ of the unit circle $\mathbb{S}^1$ through this isomorphism, $\vert\vert \cdot\vert\vert_{\nu}$ is bounded from above by some constant $C_\nu$ on $N\cap S$ (by Remark \ref{cont}, as $\mathbb{S}^1$ is compact in $(\mathbb{C}, \vert \cdot \vert)$). Henceforth, for all $(x_1,\dots, x_t)\in (N\cap S)^{\times t}$ we have 

\begin{equation}\label{bound}
\prod_{\nu\in\mathbb{V}_\infty(N)}\text{max}\{||x_1||_{\nu},\ldots ,||x_t||_{\nu}\}\leq C_K:=\prod_{\nu\in\mathbb{V}_\infty(N)}C_\nu.
\end{equation}
Since, by our hypothesis, $\frac{\sigma_{s+k}(u)}{\sigma_{s+1}(u)}\in\mathbb{S}^1$ for all $k=2,\dots, t$, we get from (\ref{bound}) and (\ref{infi}) that
$$H_N(P_u)\leq C_K$$  for all $u \in U.$ Henceforth the image of $\rho$ is finite according to Theorem \ref{NW}.

Therefore, $G_U\subset U$ is of finite index; in particular, $\text{rk}_\Z(G_U)=s$. Consider the field $k:=\{x\in K\vert \sigma_{s+1}(x)=\dots=\sigma_{s+t}(x)\}$
and denote by $(s', t')$ its signature. As at least $t$ embeddings of $K$ lie above a single embedding of $L$, we see that, according to Theorem \ref{deginf},
\begin{equation}\label{deg}
[K:k]\geq t.
\end{equation}
As all the complex embeddings of $K$ lie above a single embedding $\tau$ of $k$, and hence all real embbedings of $k$ distinct from $\tau$ lift to real embeddings of $K$, we must have
$t'\in \{0, 1\}$
and also
\begin{equation}\label{ine}
s'-1\leq \frac{s}{[K:k]}\leq \frac{s}{t}.
\end{equation}
(also by Theorem \ref{deginf}). More, from $G_U\subset \mathcal{O}_k^*$ and Dirichlet's unit theorem \cite[Chapter 5.Theorem 38]{M} we get
\begin{equation}\label{diri}
s'+t'-1\geq s.
\end{equation}

 Now, if $t'=1$, from (\ref{deg}) we get $s'+2t'\leq \frac{s}{t}+2$, that is, $s'\leq\frac{s}{t}$; but (\ref{diri}) implies $s'\geq s$ hence we get $s\leq\frac{s}{t}$ which forces $t=1$ since $s\geq 1.$ If $t'=0$, from (\ref{diri}) we get $s'-1\geq s$, so from (\ref{ine}) we get $\frac{s}{t}\geq s$ which forces again $t=1.$
\end{proof}

\begin{rem} Note that Theorem \ref{main} is also removing the admissibility hypothesis (cf. \cite{OeTo}), upon which the previous results rely. 
\end{rem}

\begin{cor}\label{qc} An Oeljeklaus--Toma manifold $X(K, U)$ is l.c.K if and only if $U$ is a congruence subgroup (see Definition \ref{cong}).
\end{cor}

\begin{proof}
 By Proposition \ref{Bat}, Theorem \ref{main}, Dirichlet's unit theorem, and Theorem \ref{chevy}.  
\end{proof}

\begin{thm}\label{pluri} Let $K$ be a number field of signature $(s, t)$ satisfying $1\leq s\leq t$ and  $U\subset\mathcal{O}_K^*$ a subgroup of rank $s$. Assume that  for any $u\in U$, $\sigma_i(u)\vert\sigma_{s+i}(u)\vert^2=1$, for all $i\in\{1,\ldots, s\}$, and $\vert\sigma_{s+i}(u)\vert^2=1$ for all $i\in\{s+1,\ldots, t\}$. Then $s=t$.
\end{thm}

\begin{proof} Suppose, on the contrary, that $t\geq s+1$. As before, consider $N=K^{\text{nor}}$ and the map
\begin{align*}
U\ni u \mapsto [\sigma_{s+i_1}(u):\ldots :\sigma_{s+i_{t-s}}(u)]\in \P^{t-s-1}(N), 
\end{align*}
where $i_1<\ldots <i_{t-s}$ are the indices from $\{s+1,\ldots, t\}$. Exactly as in the first part of the proof of Theorem \ref{main}, we get $\text{rk}_{\Z}(G_U)=s$, where now $G_U=\{u\in U; \sigma_{s+i_1}(u)=\ldots =\sigma_{s+i_{t-s}}(u)\}$. Next, consider the subextension $k=\{x\in K; \sigma_{s+i_1}(x)=\ldots =\sigma_{s+i_{t-s}}(x)\}$, and let $(s', t')=\text{sgn}(k)$.

On the other hand, by the definition of $k$, there exists an embedding $\tau$ of $k$ that has at least $t-s$ extensions to $K$.

Assume first that $\tau$ is real. Then we actually have exactly $2(t-s)$ extensions of $\tau$. Indeed, $\sigma_{s+i_j}, \overline{\sigma_{s+i_j}}$ extend $\tau$ (since $\tau$ is real), so we have at least $2(t-s)$ extensions. If there is at least one more extension, that must be either $\sigma_i$ or $\sigma_{s+i}$ (or its conjugate) for some $i\in\{1,\ldots, s\}$; if it is $\sigma_i$, then the hypothesis and construction implies that $\sigma_i(u)=\vert\sigma_i(u)\vert=1$ for all $u\in G_U$, which is absurd since this is a subgroup of rank $s\geq 1$, and if it is  $\sigma_{s+i}$, then $\vert\sigma_{s+i}(u)\vert=1$, so again $\sigma_i(u)=\sigma_i(u)\vert\sigma_{s+i}(u)\vert^2=1$ for all $u\in G_U$. We conclude that 
\begin{equation}\label{eq1}
[K:k]=2(t-s)
\end{equation}
(according to Theorem \ref{deginf}). 

Next, notice that, according to relation (\ref{eq1}), the number of possible lifts to $K$ of a complex embedding of $k$ is $2t-\nu-2(t-s)=2s-\nu$, where $\nu$ denotes the number of complex embeddings of $K$ which give a real embedding of $k$. By Theorem \ref{deginf} and relation (\ref{eq1}) there are exactly $2(t-s)$ embeddings of $K$ which extend a complex embeddings of $k$, and hence taking into account all the complex embeddings of $k$ we conclude that $2(t-s)2t'\leq2s-\nu$, i.e. $2t'\leq\frac{2s-\nu}{2(t-s)}$. For the real embeddings of $k$ distinct from $\tau$, the number of which being $s'-1$, we get $s+2t-\nu'-2(t-s)=3s-\nu'$ possible lifts, where $\nu'$ is the number of complex embeddings of $K$ which give complex embeddings of $k$. Thus, we get as before that $s'-1\leq\frac{3s-\nu'}{2(t-s)}$, hence $s'\leq\frac{s+2t-\nu'}{2(t-s)}=\frac{s+\nu}{2(t-s)}$ since, by definitions, $\nu+\nu'=2t$.

Therefore, in the case when $\tau$ is real, the relation $[K:\mathbb{Q}]=[K:k][k:\mathbb{Q}]$ gives $s+2t=2(t-s)(s'+2t')\leq s+\nu+2s-\nu=3s$, so we arrive at the contradiction $t\leq s$.  .

It remains to treat the case when $\tau$ is complex. The same reasons as above show that 
\begin{equation}\label{eq2}
[K:k]=t-s
\end{equation}
(since, this time, $\tau$ is complex, so no two conjugate emebeddings of $K$ can induce $\tau$). There are at most $2t-\nu-(t-s)=t+s$ possible lifts to $K$ of a complex embedding of $k$ distinct from $\tau$, and hence, by the relation (\ref{eq2}) and Theorem \ref{deginf}, $2t'\leq\frac{t+s-\nu}{t-s}+1=\frac{\nu'}{t-s}$. Next, a real embedding of $k$ has at most $s+2t-\nu'-2(t-s)=3s-\nu'$ extensions to $K$ because we must substract $2(t-s)$, not just $(t-s)$, since we're in the case where $\tau$ is complex, and here we are dealing with a real one. Hence, Theorem \ref{deginf} together with relation (\ref{eq2}) gives $s'\leq\frac{3s-\nu'}{t-s}$. Again, $[K:\mathbb{Q}]=[K:k][k:\mathbb{Q}]$ gives $s+2t\leq 3s$, contradicting our assumption.

Therefore, the inequality $t\geq s+1$ is impossible, and since, by hypothesis, $t\geq s$, it follows that $t=s$.
\end{proof}

\begin{rem} Though far less complicated, we notice that the previous two proofs are analogous to the proof of Mordell--Weil theorem \cite[C.0.1]{S}. Indeed, we see here that the proofs consists of two parts, one using the height in order to show that $U/G_U$ is finite ("weak Mordell--Weil") and the other using a ramification argument at infinity. 
\end{rem}

\subsection*{Acknowledgment}
We thank Vicen\c tiu Pa\c sol for a careful reading of a first draft of the paper and Victor Vuletescu for pointing us the problem.

\end{document}